\voffset-0.4in
\documentstyle{amsppt}
\magnification=1200
\document
\topmatter
\title
Holomorphic forms on threefolds
\endtitle
\author 
Tie Luo and Qi Zhang
\endauthor
\date 
\enddate
\address
Department of Mathematics, University of Texas, Arlington, TX 76019
\endaddress
\email
luo\@uta.edu
\endemail
\address
Department of Mathematics, University of Missouri, Columbia, MO 65211
\endaddress
\abstract
Two conjectures relating the Kodaira dimension of a smooth projective variety and existence of number of nowhere vanishing 1-forms on the variety are proposed and verified in dimension 3.
\endabstract
\email
qi\@math.missouri.edu
\endemail
\endtopmatter
\openup6pt
\heading
\S1 introduction
\endheading
Let $X$ be a smooth projective variety over the complex number field $\Bbb C$. Let $\omega\in H^0(X,\Omega_X^i)$ be a nontrivial global holomorphic differential $i$-form. We would like to understand the nature of zero locus of $\omega$ in terms of the birational geometry of $X$. On one hand algebraic varieties are birationally classified according to their Kodaira dimensions using high multiples of holomorphic form of top degree-the canonical bundle and Mori's theory of extremal contractions and flips are birational operations done to part of the base locus of (multiples of ) canonical bundle, on the other hand very little is known about the impact of zero locus of holomorphic forms of lower degree has on the birational geometry of the underlying variety.  The first question one may ask is: What makes $\omega$ to have zero locus? It is proved in [Z] that for any global holomorphic 1-form
$0\neq\omega\in H^0(X,\Omega_X)$, the zero locus $Z(\omega)$ is not empty provided that the canonical bundle $K_X$ is ample. It is natural to suspect that the same conclusion should hold for varieties of general type.

Indeed we propose two conjectures:
\demo{Conjecture 1} Assume $X$ is of general type and $0\neq\omega\in H^0(X,\Omega_X^1)$. Then $Z(\omega)\neq\phi$.
\enddemo
and more precisely
\demo{Conjecture 2} Assume there is a subspace $V$ of dimension $k$ in $H^0(X,\Omega_X)$ for some  $1\leq k\leq\text{dim}(X)$ and $Z(\omega)=\phi$ for any $0\neq\omega\in V$. Then the Kodaira dimension of $X$ is bounded from above by $\text{dim}(X)-k$.
\enddemo 

These conjectures are trivially true for curves. As an exercise we ask readers to check them for surfaces.
The purpose of the present paper is to verify that they also hold in dimension $3$ by proving: 

\proclaim{Theorem 1} Let $X$ be a threefold of general type. For any global holomorphic 1-form
$0\neq\omega\in H^0(X,\Omega_X)$, its zero locus $Z(\omega)$ is not empty.   
\endproclaim

and 

\proclaim{Theorem 2} Let $X$ be a threefold admitting a $k$-dimensional subspace $V\subset H^0(X,\Omega_X^1)$ such that $z(\omega)=\emptyset$ for any $0\neq\omega\in V$. Then the Kodaira dimension of $X$ is at most $3-k$.
\endproclaim

Some of the techniques used in proving these results can not be generalized to higher dimensions. We plan to discuss similar results in higher dimensions using different tools in a future paper.

In contrast to the result on 1-forms, there are examples (see section 4 below) of threefolds with ample canonical bundle and carrying nowhere vanishing 2-forms. Using a completely different approach and viewing a 2-form as a meromorphic vector field, we have

\proclaim{Theorem 3} Let $X$ be a threefold of general type. $\omega\in H^0(X,K_X\otimes\Omega^2_X)$ is a canonically twisted holomorphic 2-form. Then the zero locus $Z(\omega)$ of $\omega$ is not empty.
\endproclaim

\demo{Acknowledgments}. We thank S. Keel, J. McKernan, and M. Reid for their interest in and comments on our results. 
\enddemo

\heading
\S 2 Proof of Theorem 1
\endheading

The following result uses the classification of extremal contractions in dimension 3. It is related to the ideas used in [L] and reduces the proof of Theorem 1 to those $X$ admitting smooth minimal models.

\proclaim{Lemma 2.1} Assume $X$ is smooth and of  general type. If $X$ does not admit a smooth minimal model, then every holomorphic 1-form on $X$ has nonempty zero locus.
\endproclaim
\demo{Proof} Since $K_X$ is not nef, there is an extremal contraction $\pi: X\rightarrow Y$. Let $E$ be the exceptional locus. Mori [M] says that $f$ is divisorial with the exceptional divisor $\Bbb P^2$ with normal bundle $\Cal O(-1)$ or $\Cal O(-2)$ , $\Bbb P^1\times\Bbb P^1$ with normal bundle $\Cal O(-1)\otimes\Cal O(-1)$, a cone $E$ over rational normal curve of degree 2 with normal bundle $\Cal O(-1)$, a $\Bbb P^1$ bundle over a smooth curve.

One has the following exact sequence
$$
0\rightarrow N_{E/X}^*\rightarrow \Omega_X\rightarrow \Omega_E\rightarrow 0.
$$

Note that unless $E$ is a  $\Bbb P^1$ bundle over a smooth curve, $H^0(E,\Omega_E)=0$. This implies that restrictions of a nontrivial global 1-form on $E$ provides a nontrivial global section of positive line bundle $N_{E/X}^*$, which has nontrivial zero locus.

When $E$ is a  $\Bbb P^1$ bundle over a smooth curve, the resulting $Y$ is smooth. From our assumption, $K_Y$ is not nef. One can repeat the argument done for $X$ and the process has to terminate after finitely many steps. So every holomorphic 1-form on $X$ has nonempty zero locus.

\qed
\enddemo
Remark. If $X$ admits a 1-form with isolated zeros, $X$ must be minimal. In [L] we showed that it is possible to have a nonminimal threefold carrying a 2-form with isolated zeros.

Thus we may assume that the canonical bundle of $X$ is nef and big. Let $f: X\rightarrow Y$ be the birational morphism to the canonical model of $X$. $Y$ is known to have Gorenstein canonical singularities. Reid [R1] constructed explicitly a crepant resolution for such singularities. Let $g:X'\rightarrow Y$ be the resolution such that $g$ is crepant and $X'$ has at most isolated cDV points. Since both $K_X$ and $K_{X'}$ are nef, it is known that $X$ and $X'$ are different by flops. One way to prove Theorem 1 is to show that any 1-form on $X'$ has to have nontrivial zero locus and flops do not change the existence of the zero locus.

Here we take a different approach using Kodaira-Akizuki-Nakano type vanishing theorem for semismall morphisms proved by Migliorini as in [Mi]. A birational morphism $\phi$ is called \it semismall \rm if for every irreducible component $E$ in the exceptional locus, 
$$
\text{codim}E\geq \text{dim}\phi^{-1}(p)
$$
for any $p\in\phi (E)$. For example in dimension 3 a semismall morphism contracts curves to points or divisor to curves.

\proclaim{Proposition 2.2(Theorem 4.6 in [Mi])} Let $X$ be a smooth projective variety over $\Bbb C$ with a semiample divisor $L$. Assume the morphism $\phi_{|mL|}$ is semismall for some $m>>0$. Then
$$
H^p(X, L^{-1}\otimes \Omega_X^q)=0
$$
for $p+q<$ dim$X$.

\endproclaim

The following general result deals with maps which contract divisors to points.

\proclaim{Lemma 2.3} Let $\pi: X\rightarrow Y$ be a projective birational morphism with $X$ smooth and $Y$ normal. Assume an irreducible divisor $E\subset X$ is contracted to a point. Let $\omega\in H^0(X,\Omega_X)$. Then $Z(\omega)$ is not trivial along $E$.
\endproclaim
\demo{Proof} Assume $\pi(E)=p$. Let us embed a neighborhood of $p$ in some $\Bbb C^N$ so that $p$ is mapped to the origin $O$ and still call the composition $\pi$. For any $q\in E$, one has $\pi^*: m_{O,\Bbb C^N}/m_{O,\Bbb C^N}^2 \rightarrow
m_{q,X}/m_{q,X}^2$ between the differentials.

If $q$ is a singularity of $E$ and assume $E$ is locally defined by some $f$, then $f\in m_{q,X}^2$. This implies $\pi^*=0$ and $q\in z(\omega)$. If $E$ is smooth at $q$, then the rank of $\pi^*$ is either $0$ or $1$. When the rank is 0, again $q\in z(\omega)$. So we may reduce the proof of the desired result to the case in which the rank of $\pi^*$ is 1 along $E$ (in particular $E$ is smooth and $\pi^*$ of any element of $m_{O,\Bbb C^N}/m_{O,\Bbb C^N}^2$ is a scalar multiple of the class $x_1+ m_{q,X}^2$ where $x_1$ locally defines $E$. The map $ m_{q,X}/m_{q,X}^2 \rightarrow
m_{q,E}/m_{q,E}^2$ is simply dropping $x_1$ component, hence sending a scalar multiple of the class $x_1+ m_{Q,X}^2$ to 0).

In this case $q\notin Z(\omega)$ for any $q\in E$ implies that $\omega$ defines a nowhere vanishing section of $N_{E/X}^*$. This is against the adjuction lemma of Kawamata in [K].

\qed
\enddemo 

Now we can finish the proof of Theorem 1.

\proclaim{Theorem 2.4=Theorem 1} Let $X$ be a threefold of general type. For any global holomorphic 1-form
$0\neq\omega\in H^0(X,\Omega_X)$, its zero locus $Z(\omega)$ is not empty.   
\endproclaim
\demo{Proof} If $K_X$ is not nef, Lemma 2.1 says $Z(\omega)$ is not empty. Assume $K_X$ is nef, then $K_X$ is semiample. Let $\phi_m=\phi_{mK_X}$ for $m>>0$. By Lemma 2.3 we may assume $\phi_m$ is semismall. Assuming $Z(\omega)$ is empty, one has exact sequence

$$
0\rightarrow \Cal O\overset\omega\to\rightarrow \Omega^1_X\overset\wedge\omega\to\rightarrow\Omega_X^2\overset\wedge\omega\to\rightarrow\Omega_X^3\rightarrow 0.
$$

Arguing as in [Z], using Proposition 2.2 when $L=K_X$ one gets an embedding of $\Omega_X^3$ into $\Omega_X^2$. It follows that $H^0(X,K_X^{-1}\otimes\Omega_X^2)\neq 0$ which is impossible again by 2.2.

\qed
\enddemo
It is easy to construct a threefold of general type with ample canonical bundle admitting a 1-form whose zero locus is of dimension 0, 1, or 2.

\heading
\S 3 Proof of Theorem 2
\endheading

On a threefold $X$, the maximal dimension of a subspace $V\subset H^0(X,\Omega_X^1)$ with the property that any $0\neq\omega\in V$ has no zero locus is 3 because the rank of $\Omega_X^1$ is $3$. By taking a product of a variety of general type and an Abelian variety, one gets examples of threefold of Kodaira dimension $3-k$ and having a $k$-dimensional subspace of 1-forms with the property that any nonzero element has no zero.

\proclaim{Theorem 3.1=Theorem 2} Let $X$ be a threefold admitting a $k$-dimensional subspace $V\subset H^0(X,\Omega_X^1)$ such that $z(\omega)=\emptyset$ for any $0\neq\omega\in V$. Then the Kodaira dimension of $X$ is at most $3-k$.
\endproclaim
\demo{Proof}
When $k=1$, the claim is implied by Theorem 1. When $k=3$, let $\omega_1, \omega_2, \omega_3$ be a basis of $V$, $\omega_1\wedge\omega_2\wedge\omega_3$ provides a nowhere zero section of $K_X$. So $K_X$ is trivial and $\kappa (X)=0$.

Let us assume $k=2$. By Lemma 2.1 we see that $X$ is minimal because of existence of $V$. We also know $\kappa (X)\neq 3$ by Theorem 1. Thus it is enough to show $\kappa (X)\neq 2$. Assuming the contrary holds. Let $\phi:X\rightarrow S$ be the IItaka fibration defined by high multiples of $K_X$. $\phi^{-1}(s)$ is an elliptic curve for a general $s\in S$. $S$ is surface of general type with at worst canonical singularities. Let $\alpha: X\rightarrow \text{Alb}(X)$ be the Albanese map. The existence of $V$ implies that $\alpha(X)$ is at least of dimension 2.

We consider two cases:

Case I. $\alpha (\phi^{-1}(s))$ is a point for a general $s\in S$. In this case $\alpha (X)$ has dimension 2. Let $X\rightarrow Y\rightarrow \alpha (X)$ be the Stein factorization of $\alpha$. It is easy to see that $S$ and $Y$ are birationally equivalent. In particular 1-forms on $X$ are pullbacks of 1-forms (on the smooth locus) of $S$. Assume the subspace $V$ has a basis $\omega_1, \omega_2$, which correspond to $\omega_1', \omega_2'$ on $S$. $\omega_1', \omega_2'$ span a 2-dimensional subspace $V'$ of $H^0(S,\Omega_S^1)$ having the property that for any $0\neq\omega'\in V'$ the zero locus on the smooth locus of $S$ is empty. Then $\omega_1'\wedge\omega_2'$ generates a nowhere zero section of $K_S$ on the smooth locus. Hence $\kappa(S)=0$, a contradiction.

Case II. $\alpha (\phi^{-1}(s))$ is a curve for a general $s\in S$. In this case $0=\kappa(\phi^{-1}(s))\geq \kappa(\alpha (\phi^{-1}(s))\geq 0$. So $\alpha (\phi^{-1}(s))=E_s$ is also an elliptic curve in $A=\text{Alb}(X)$. Since $A$ contains at most countably many sub Abelian varieties of dimension 1, $E_s$ are translations of a particular sub Abelian variety $A_1$ of dimension 1 . Let $\pi: A\rightarrow A/A_1$ be the quotient map. One has the following diagram:
$$
\CD
X@>\alpha>> A\supset \alpha (X)\\
@V\phi VV @VV\pi V\\
S@. A/A_1\supset \pi(\alpha(X))
\endCD
$$
If $\alpha (X)$ is of dimension 3, $\pi(\alpha (X))$ is of dimension 2. Again let $X\rightarrow Y\rightarrow \pi(\alpha (X))$ be the Stein factorization of $\pi\circ\alpha$. $S$ is birationally equivalent to $Y$ and we are done by case I. If $\alpha (X)$ is of dimension 2, $C=\pi(\alpha (X))$ is of dimension 1 in the Abelian variety $A/A_1$. Let $\omega_1', \omega_2'$ be the 1-forms on $A$ which pullback to a basis of $V$. There is a nontrivial linear combination $\omega'=c_1\omega_1'+c_2\omega_2'$ which is identically zero on the fibers of $\pi$. $\omega'$ is a pullback of some $\omega''$ on $A/A_1$. Assume $C$ is of general type and we may assume $C$ is smooth, $\omega"|C$ vanishes at some point $p\in C$, which in turn implies $c_1\omega_1+c_2\omega_2$ vanishes on the fiber of $\pi\circ \alpha$ in $X$ over $p$. If $C$ is smooth elliptic curve, the pullback of $\omega'=c_1\omega_1'+c_2\omega_2'$ (which is the same as the pullback of $\omega''$) would have a zero unless the morphism $\pi\circ\alpha$ is smooth and $\omega''$ does not vanish on $C$. So let us assume $\pi\circ\alpha$ is smooth. Note that each fiber the morphism has to be a surface of general type in order to keeping $\kappa (X)=2$. Furthermore each fiber has to be minimal since $K_X$ is nef. Thus $\pi\circ\alpha: X\rightarrow C$ represents a smooth family of minimal surfaces of general type over an elliptic curve. It has to be isotrivial thanks to Theorem 1.1 in [Mi]. So up to an \'etale covering, X is isomorphic to the product of $C$ and a fiber of $\pi\circ\alpha$. It implies that $V$ contains a nonzero member $\omega$ coming from a 1-form on the fiber which is a surface of general type. As any 1-form on a surface of general type has to have a zero, $Z(\omega)$ contains a curve. A contradiction.

\qed
\enddemo

\heading
\S 4 2-forms on threefolds of general type
\endheading

After the discussion regarding zeros of 1-forms, it is quite natural to wonder whether similar results on 1-forms still hold true for higher orders forms, i.e;

Question: Let $X$ be a threefold of general type. Is it true that any global differential 2-form on $X$ must have nonempty zero locus?

Campana and Peternell [CP] constructed the following example which shows that a global section of $\Omega_X^2$ may not have zero locus even when $K_X$ is ample:

\demo{Example} Let $Y$ be projective holomorphic symplectic fourfolds. Let $X$ be a hypersurface of sufficiently high degree. Then the nondegenerate 2-form on $Y$ provides a no where vanishing 2-form on $X$.
\enddemo

However by virtue of [L] it is known that when $K_X$ is not nef, any 2-form has nontrivial zero locus. It would be interesting to see whether there are examples with nef canonical bundle and a 2-form without any zero. 

Here we will show that any canonically twisted 2-form $\omega\in H^0(X,  K_X\otimes\Omega_X^2)$ on a threefold of general type has nontrivial zero locus.

Let $X$ be a smooth projective variety of dimension $n$. $\theta\in H^0(X,T_X\otimes L)$ with isolated zeros. Then locally around each zero $p$, under coordinates $z_1, z_2,...,z_n$, one can express $\theta$ as
$$
\theta=\sum_{i=1}^nf_is\otimes\frac{\partial}{\partial z_i}
$$
where $f_i$ are holomorphic functions and s is a local section of $L$. The Jacobian at the zero is defined as 
$$
A_p=(\frac{\partial f_i}{\partial z_j})_{|p}
$$
which is an $n\times n$ matrix. $A_p$ is invariant under coordinate change upto a nowhere zero holomorphic function. $p$ is nondegenerate if $det(A_p)\neq 0$.

The following lemma is a corollary of the residue theorem due to Bott (see [C]).

\proclaim{Lemma 4.1} Notations as above and assume $\theta$ has only nondegenerate zeros. Then

$$
\int_X\bar c_1(T_X\otimes L)^n=\sum_{\text{$p$ a zero}} \frac{(Tr(A_p))^n}{det(A_p)} 
$$
and
$$
\int_X\bar c_n(T_X\otimes L)=\sum_{\text{$p$ a zero}}\frac{det(A_p)}{det(A_p)}
$$

where $\bar c_r(T_X\otimes L)=c_r(T_X)+c_{r-1}c_1(L)+...+c_1(L)^r$.
\endproclaim

Remark. Indeed Chern's proof [C] shows that away from the zero locus, $\bar c_1(T_X\otimes L)^n$ and $\bar c_n(T_X\otimes L)$ can be written as derived differential forms. This fact is used in the proof of the following 

\proclaim{Theorem 4.2=Theorem 3} Let $X$ be a threefold of general type. $\omega\in H^0(X, K_X \otimes \Omega^2_X)$ is a canonically twisted holomorphic 2-form. Then the zero locus $Z(\omega)$ of $\omega$ is not empty.
\endproclaim
\demo{Proof} We may assume $K_X$ is nef since otherwise there is a divisorial extremal contraction and an analysis similar to that in section 2 of [L] shows that for all possible divisorial extremal contraction any canonically twisted 2-form has nontrivial zero locus along the exceptional divisor.

Since $\Omega^2_X$ is isomorphic to $K_X\otimes T_X$, $\omega$ corresponds to a
 global section $\theta$  in $H^0(X,2K_X\otimes T_X)$, which can be viewed as a meromorphic vector field
with coefficient in the line bundle $2K_X$. Assume $\theta$ (or equivalently $\omega$) does not have zero, the remark about the above lemma says
$$
\bar c_1(2K_X\otimes T_X)^3=0.
$$
However the left hand side is
$$
\bar c_1(2K_X\otimes T_X)^3=(c_1(T_X)+c_1(2K_X))^3=K_X^3>0,
$$
a contradiction.

\qed
\enddemo

Finally the geometric meaning of $c_3(\Omega^2_X)$ (which is $c_3-c_1c_2$ of $X$) in terms of singularities appearing on a minimal model is explained in:
\proclaim{Theorem 4.3} Let $X$ be a threefold of general type. Assume $\Omega^2_X$ admits a section $\omega$ with isolated nondegenerate zeros. Then $c_3(\Omega^2_X)$ is greater or equal to the number of singularities whose indices are bigger than $1$ on a minimal model.
\endproclaim
\demo{Proof}
Let $X_{min}$ be a minimal model of $X$. Assuming $\omega$ is a 2-form with only isolated nondegenerate zeros on $X$. By [L], $\omega$ corresponds to a (global locally invariant) 2-form on $X_{min}$ not vanishing at the singularities whose indices are bigger than $1$. Moreover every such singularity produces at least one isolated zero on $X$. Bott's residue formula, when taking the top Chern class, says
$$
\bar c_3(K_X\otimes T_X)=\sum_{\text{$p$ a zero}}\frac{det(A_p)}{det(A_p)}=\text{number of zeros of $\omega$}.
$$
This implies
$$
\bar c_3(K_X\otimes T_X)=c_3(T_X)+c_2(T_X)c_1(K_X)+c_1(T_X)c_1(K)^2+c_1(K_X)^3=c_3-c_1c_2
$$
$$
\geq \text{number of singularities whose indices are bigger than 1}.
$$
\qed
\enddemo


\Refs

\ref 
\by [C] Chern,~S.S
\paper
Meromorphic vector fields and characteristic numbers 
\jour Scripta Mathematica
\vol 29
\pages 243-251
\yr 1973
\endref

\ref
\key 
\by [CP] Campana,~F and Peternell,~T
\paper Holomorphic 2-forms on complex threefolds 
\yr 2000
\vol 9
\pages 223-264
\jour Journal of Algebraic Geometry
\endref


\ref 
\by [K] Kawamata,~Y
\paper
On the length of an extremal rational curve
\jour Invent. Math.
\pages 609-611
\vol 105
\yr 1991
\endref

\ref 
\by [L] Luo,~T
\paper
Global 2-forms and pluricanonical systems on threefolds
\jour Math. Ann.
\pages to appear
\yr 2000
\endref 

\ref
\by [Mi] Migliorini, ~L, A smooth family of minimal surfaces of general type over a curve of genus at most one is trivial
\jour J. Algebraic Geom.
\vol 4
\yr 1995
\pages 356-363
\endref

\ref
\no 
\by [Mo] Mori,~S
\pages 133-176
\paper Threefolds whose canonical bundles are not numerically
effective
 \yr 1982 \vol 116
\jour Ann. of Math.
\endref
  
\ref 
\by [R1] Reid,~M 
\paper
Canonical 3-folds
\inbook Journ\'ee de G\'eometri alg\'ebrique d'Angers
\pages 273-310
\vol 
\yr 1980
\endref 

\ref 
\by [R2] Reid,~M 
\paper
Minimal models of canonical 3-folds
\inbook Algebraic Varieties and Analytic Varieties, Adv. Stud. Pure Math.
\pages 131-180 
\vol 1
\yr 1980
\endref 

\ref 
\by [R3] Reid,~M 
\paper
Young person's guide to canonical singularities
\inbook Algebraic Geometry-Bowdoin 1985
\pages 333-414
\vol 46
\yr 1987
\endref 




\ref 
\by [U] Ueno,~K 
\paper
Classification theory of algebraic varieties and compact complex spaces
\vol 439
\jour Lecture Notes in Mathematics
\yr 
\endref    

\ref 
\by [Z] Zhang,~Q 
\paper
Global holomorphic one-forms on projective manifolds with ample canonical bundles
\jour Journal of Algebraic Geometry
\pages 777-787
\vol 6
\yr 1997
\endref 

 








\endRefs
\enddocument